\newfont{\Bbb}{msbm10 scaled\magstep1}

\newcommand{\ZZ}{{\mathbb Z}}
\newcommand{\CC}{{\mathbb C}}
\newcommand{\QQ}{{\mathbb Q}}
\newcommand{\NN}{{\mathbb N}}

\def\TT{{\bf T}}

\def\ulamb{{\underline{\lambda}}}
\def\umu{{\underline{\mu}}}

\def\lra{\longrightarrow}

\def\ee{\epsilon}
\def\proof{{\bf Proof.}\,\,}
\def\qed{\hfill{{\bf QED}}}
\def\qed{\hfill\vrule height4pt width4pt depth0pt}  

\def\be{\begin{equation}}
\def\ee{\end{equation}}
\def\ep{{\epsilon}}
\def\kformep{\epsilon^1\wedge\ldots\wedge\epsilon^k}

\def\ikformep{\epsilon^{i_1}\wedge\ldots\wedge\epsilon^{i_k}}

\def\ihikformep{\epsilon^{i_1+h_1}\wedge\ldots\wedge\epsilon^{i_k+h_k}}
\def\rikformep{\epsilon^{1+r_1}\wedge\ldots\wedge\epsilon^{k+r_k}}

\def\wM{\bigwedge M}

\def\wkM{\bigwedge^k M}

\documentclass[12]{article}

\usepackage{amsfonts,amssymb,theorem}

\usepackage{graphicx}
\usepackage{eufrak}

\newfont{\eulerfraktur}{eufm10 scaled\magstep1}

\theoremstyle{change}
\theorembodyfont{\rmfamily}
\global
\global

\newtheorem{thm}{Theorem}[section]

\newtheorem{prop}{Proposition}[section]
\newtheorem{corol}{Corollary}[section]
\newtheorem{lemma}{Lemma}[section]

\newtheorem{claim}{}[section]

\begin{document}
\title{\normalsize{\Large Schubert Calculus via Hasse-Schmidt Derivations}
}

\author{{\small LETTERIO GATTO}\thanks{Work partially
sponsored by MIUR  (Progetto Nazionale ``Geometria sulle Variet\`a
Algebriche",  coordinatore Sandro Verra), and supported by GNSAGA-INDAM.}\\Ê
\\ {Dipartimento di Matematica,
Politecnico di Torino,}
\\
{\normalsize Corso Duca degli Abruzzi 24, 10129 Torino -- (ITALY)}\\
}
\date{}

\maketitle
\begin{abstract}
\noindent
A natural Hasse-Schmidt derivation on the exterior algebra of a free module realizes the (small quantum) cohomology ring of the grassmannian $G_k(\CC^n)$ as a ring of operators on the exterior algebra of a free module of rank $n$. Classical Pieri's formula can be interpreted as Leibniz's rule enjoyed by special Schubert cycles with respect to the wedge product.
\end{abstract}

\section{Introduction}

The main purpose of this note is to suggest a new simple point of view to look at (small quantum) Schubert Calculus, based on elementary considerations of linear algebra.
To get into the matter of the paper, it seems worth to start with an example.
Let $D$ be the endomorphism of $M_4:=\oplus_{1\leq i\leq 4}\ZZ\cdot\epsilon^i$ defined by 
$D\epsilon^i=\epsilon^{i+1}$, for $1\leq i<4$, and $D\epsilon^4=0$. Extend it to $\bigwedge^2M_4$, by imposing Leibniz's rule with respect to $\wedge$, and  compute $D^4(\epsilon^1\wedge\epsilon^2)$. One has:
\begin{eqnarray*}
&{ }&D^4(\epsilon^1\wedge\epsilon^2)=D\circ D\circ D\circ D(\epsilon^1\wedge\epsilon^2)=D\circ D\circ D(\epsilon^1\wedge\epsilon^3)=\\&{Ê}&
D\circ D(\epsilon^2\wedge\epsilon^3+\epsilon^1\wedge\epsilon^4)=D(2\epsilon^2\wedge \epsilon^4)=2D(\epsilon^2\wedge\epsilon^4)=2\cdot\epsilon^3\wedge\epsilon^4.
\end{eqnarray*}
The claim is that the above iteration of $D$ computes the number ($=2$) of lines  intersecting four others in general position in the projective $3$-space (see e.g.~\cite{KlLa}, p.~1068--1069, 1073--1074, \cite{GH}, p.~206).
The  reason is that {\em the cohomology  ring of the grassmannian $G_k(\CC^n)$ can be realized as a natural commuative ring of endomorphisms of the $k$-th exterior power of a free module of rank $n$} (Theorem~\ref{thm3.2}). This is a consequence of the following nicer and more general fact. Let $M$ be a free $\ZZ$-module. Using a terminology borrowed from commutative algebra, as e.g. in~\cite{Mats}, p. 207, one says that $D_t:=\sum_{i\geq 0}D_it^i:\wM\lra(\wM)[[t]]$ $(D_i\in End_{\ZZ}(\wM))$ is a {\em Hasse-Schmidt derivation} on $\wM$ if it is a $\ZZ$-algebra homomorphism, i.e. if:
\be
D_t(\alpha\wedge\beta)=D_t(\alpha)\wedge D_t(\beta),\qquad \forall\alpha,\beta\in\wM\label{eq:intr01}.
\ee
Let ${\cal E}:=(\epsilon^1,\epsilon^2,\ldots)$ be a (countable infinite) $\ZZ$-basis of  a free $\ZZ$-module $M$. If $D_t$ is the unique $HS$-derivation on $\wM$ such that $D_t(\epsilon^j)=\sum_{i\geq 0}\epsilon^{i+j}t^i$ (thinking of $M$ as a submodule of $\wM$), then {\em  Schubert Calculus of $G_k(\CC^n)$, for all $(k, n)$ at once ($0\leq k\leq n$), is a formal consequence of formula~{\rm (\ref{eq:intr01})}}. This is why $D_t$ is named {\em Schubert derivation} (Def.~\ref{sdef}).

Indeed, for all $k\geq 0$,  $\wkM$ is a $D_h$-invariant submodule of $\wM$, for each  ``coefficient" $D_h$ of $D_t$; the point is that the entries of the (infinite) matrix of ${D_h}_{|_{\wedge^kM}}$ with respect to the basis $\{\ikformep: 1\leq i_1<i_2<\ldots<i_k\}$ of $\wkM$, can be computed via {\em Pieri's formula for ${\cal S}$-derivations} (Theorem~\ref{Pieri}):
\[
D_h(\ikformep)=\sum\ihikformep
\] 
the sum being over all non-negative $(h_1,\ldots,h_k)$ such that $h_1+\ldots+h_k=h$ and
\[
1\leq i_1\leq i_1+h_1<i_2\leq i_2+h_2<\ldots<i_{k-1}\leq i_{k-1}+h_{k-1}< i_k .
\]
This is precisely classical Pieri's formula, as briefly explained in Sect.~\ref{concrem}.
 
 
Let $M_n$ be the submodule of $M$ spanned by $(\epsilon^1,\ldots,\epsilon^n)$.  Via the formal identification $\rikformep\mapsto\sigma_\ulamb\cap [G_k(\CC^n)]$ (the Schubert cycle $\sigma_\ulamb$ corresponding to the partition $\ulamb=(r_k,\ldots,r_1)$ capped with the fundamental class of the Grassmannian) and using the Chow basis theorem for the cohomology of $G_k(\CC^n)$, one concludes that, in fact, {\em the cohomology ring of $G_k(\CC^n)$ is a (commutative) ring of endomorphisms on $\wkM_n$} and that all such, varying $k$ and $n$, are quotient of a (same) natural ring of derivations on $\wM$ (Thm.~\ref{thm3.2}). 

The results of this work have been recently improved and generalized by Laksov and Thorup~(\cite{LakTh}) to grassmannian bundles, using the theory of symmetric functions and of splitting algebras, allowing them to study, in general, the cohomology of (partial) flag varieties of a finite dimensional vector space over an algebraically closed field  (\cite{LakTh2}).

The current shape of this paper is mostly due to the patients corrections and substantial remarks of the Referees whom the author is very grateful.  
For warm encouragement, but  especially for his as friendly as sharp criticism, no word would be enough to thank Dan Laksov.

\section{Schubert DerivationsÊ}\label{Sect3}

Let $\wM=\bigoplus_{k\geq 0}\wkM$ be the exterior algebra of a $\ZZ$-module $M$ freely generated by ${\cal E}=(\epsilon^1,\epsilon^2,\ldots)$. Denote by 
\[
\wedge^k{\cal E}:=\{(\ikformep):1\leq i_1<i_2<\ldots<i_k\}
\] 
the induced basis of $\wkM$.

\begin{claim}\label{sdef}{\bf Definition.} {\em
A {\em Hasse-Schmidt ($HS$) derivation} on $\wM$ is a $\ZZ$-algebra homomorphism $D_t:=\sum_{i\geq 0}D_it^i:\wM\lra(\wM)[[t]]$ ($D_i\in End_{\ZZ}(\wM)$).
}
\end{claim}
Formally, the $\ZZ$-algebra homomorphism condition reads as:
\be
D_t(\alpha\wedge\beta)=D_t(\alpha)\wedge D_t(\beta),\qquad \forall \alpha,\beta\in \wM.\label{eq:HS01}
\ee
Clearly, $D_t$ is uniquely determined by its values on the elements of the basis ${\cal E}$ of $M$ (thought of as a submodule of $\wM$). 
Let $D:=(D_0,D_1,\ldots)$ be the sequence of {\em coefficients} of $D_t$. 
 Formula~(\ref{eq:HS01}) can be then rephrased by saying that each $D_h$ satisfies Leibniz's rule for $h$-th order derivatives:
\be
D_h(\alpha\wedge\beta)=\sum_{\matrix{h_1+h_2=h\cr h_i\geq 0}}D_{h_1}\alpha\wedge D_{h_2}\beta.\label{eq:LR}
\ee
In fact, the r.h.s of~(\ref{eq:LR}) is precisely the coefficient of $t^h$ in the expansion of the r.h.s. of~(\ref{eq:HS01}).

\begin{claim}\label{defscder}{\bf Definition.} {\em
The (${\cal E}$)--{\em Schubert derivation} (${\cal S}$-derivation) is the unique $HS$-derivation on $\wM$ such that
\be
D_t(\epsilon^i)=\sum_{j\geq 0}\epsilon^{i+j}t^j.\label{eq:for1}
\ee
}
\end{claim}
Such a ${\cal S}$-derivation exists: it suffices to extend a map $D_t:M\lra M[[t]]$ satisfying~(\ref{eq:for1})  to all $\wM$ by imposing~(\ref{eq:HS01}).

\noindent
Next task is to find the components of the endomorphisms $D_h:\wM\lra\wM$ ($h\geq 1$) with respect to the basis $\bigwedge{\cal E}=\cup_{k\geq 0}\wedge^k{\cal E}$. One first puts~(\ref{eq:LR}) in a more explicit form.
\begin{claim}{\bf Proposition.}\label{propprepie} {\em
For each $h\geq 0$ and  each $k\geq 1$, one has:
\be
D_h(\epsilon^{i_1}\wedge \epsilon^{i_2}\wedge\ldots\wedge \epsilon^{i_k})=
\sum_{{\matrix{_{h_1+\ldots+h_k=h}\cr _{h_i\geq 0}}}} \epsilon^{i_1+h_1}\wedge
\epsilon^{i_2+h_2}
\wedge\ldots\wedge \epsilon^{i_k+h_k}.\label{eq:for002}
\ee
}
\end{claim}
\proof
For $k=1$, formula~(\ref{eq:for002}) is Definition~\ref{defscder}. Assume it holds for $k-1$. Application of~(\ref{eq:LR}) gives:
\be
D_h( \epsilon^{i_1}\wedge \epsilon^{i_2}\wedge\ldots\wedge \epsilon^{i_k})=
\sum_{h_1=0}^{h} \epsilon^{i_1+h_1}\wedge
D_{h-h_1}( \epsilon^{i_2}\wedge\ldots\wedge \epsilon^{i_k}),
\label{eq:prepier}
\ee
where
\[
D_{h-h_1}( \epsilon^{i_2}\wedge\ldots\wedge \epsilon^{i_k})=
\sum_{h_2+\ldots+h_k=h-h_1} \epsilon^{i_2+h_2}\wedge\ldots\wedge
 \epsilon^{i_k+h_k},
\]
by the inductive hypothesis.
Thus, the right hand side of formula~(\ref{eq:prepier}) turns into:
\[
\,\,\,\,\,\,\,\,\,\,\,\,\,\,\,\, \,\,\,\,\,\,\,\,\,\,D_h( \epsilon^{i_1}\wedge \epsilon^{i_2}\wedge\ldots\wedge \epsilon^{i_k})=
\sum_{h_1+\ldots+h_k=h} \epsilon^{i_1+h_1}\wedge\ldots
\wedge \epsilon^{i_k+h_k}.\,\,\,\,\,\,\,\,\,\,\,\,\,\,\,\,\,\,\,\,\,\,\,\,\,\,\qed
\]

\noindent
Proposition~\ref{propprepie} clearly implies that $D_iD_j=D_jD_i$ for all $i,j\geq 0$.
Hence the evaluation morphism $E_D:\ZZ[\TT]\lra End_{\ZZ}(\wM)$, gotten by sending $T_i\mapsto D_i$ is well defined and maps onto the commutative subalgebra $\ZZ[D]\subset End_{\ZZ}(\wM)$ generated by $D:=(D_1,D_2,\ldots)$.  Indeed, for each $k\geq 1$, $\ZZ[D]$ can be seen as a subalgebra of $End_{\ZZ}(\wkM)$, because Definition~\ref{sdef} and/or Proposition~\ref{propprepie} imply that $D_n(\wkM)\subseteq \wkM$,  for each $n\geq 0$. 

\begin{claim}{\bf Theorem.}\label{Pieri}
{\em Let $I:=(1\leq i_1<i_2\ldots<i_k)$ be a sequence of integers. Then {\em Pieri's formula} for ${\cal S}$-derivations holds:
\be
D_h(\ikformep)=\sum_{\matrix{_{(h_i)\in H(I,h)}}} \ihikformep,\label{eq:pieruno}
\ee
where, to shorten notation, one denotes by $H(I,h)$ the set of all $k$-tuples $(h_i)$ of non-negative integers such that
\be
1\leq i_1\leq i_1+h_1<i_2\leq\ldots \leq i_{k-1}+h_{k-1}< i_k.\label{eq:pierdineq}
\ee
and $h_1+\ldots+h_k=h$
}
\end{claim}
\proof
By induction on the integer $k$. For $k=1$,  formula~(\ref{eq:pieruno})  is trivially true.
Let us prove it directly for $k=2$. For each $h\geq 0$, let us split sum~(\ref{eq:prepier}) as:
\begin{eqnarray}
D_h(\epsilon^{i_1}\wedge\epsilon^{i_2})=\sum_{h_1+h_2=h}\epsilon^{i_1+h_1}\wedge\ep^{i_2+h_2}
={\cal P}+\overline{\cal P}.\label{eq:piercorr1}
\end{eqnarray}
 where 
 \[
 {\cal P}=\sum_{\matrix{_{i_1+h_1< i_2}\cr _{h_1+h_2=h}}}\ep^{i_1+h_1}\wedge\ep^{i_2+h_2} \quad {\rm and}\quad  
 \overline{\cal P}=  \sum_{\matrix{_{i_1+h_1\geq i_2}\cr _{h_1+h_2=h}}}\ep^{i_1+h_1}\wedge\ep^{i_2+h_2}.
 \]

One contends that $\overline{\cal P}$ vanishes. In fact, on the finite set of all integers $i_2-i_1\leq a\leq i_2-i_1+h$, define the bijection $\rho(a)=i_2-i_1+h-a$. 
Then:
\begin{eqnarray*}
&{}&2\overline{\cal P}=\sum_{h_1=i_2-i_1}^h\ep^{i_1+h_1}\wedge\ep^{i_2+h-h_1}+
\sum_{h_1=i_2-i_1}^h\ep^{i_1+\rho(h_1)}\wedge\ep^{i_2+h-\rho(h_1)}=\\ &=&\sum_{h_1=i_2-i_1}^h\ep^{i_2+h-h_1}\wedge\ep^{i_1+h_1}-\sum_{h_1=i_2-i_1}^h\ep^{i_1+h_1}\wedge\ep^{i_2+h_2}=0,
\end{eqnarray*}
hence $\overline{\cal P}=0$ and~(\ref{eq:pieruno}) holds for $k=2$.  Suppose now that~(\ref{eq:pieruno}) holds for all $1\leq k'\leq k-1$. Then, for each $h\geq 0$:

\[
D_h(\ikformep)=\sum_{{h'_k}+h_k=h}D_{h'_k}(\ep^{i_1}\wedge\ldots\wedge\ep^{i_{k-1}})\wedge D_{h_k}\ep^{i_k},
\]
and, by the inductive hypothesis:
\be\sum_{(h_i)}(\ep^{i_1+h_1}\wedge\ldots\wedge\ep^{i_{k-2}+h_{k-2}}\wedge \ep^{i_{k-1}+h_{k-1}})\wedge \ep^{i_{k}+h_{k}},\label{eq:piercorr2}
\ee
summed over all $(h_i)$ such that $h_1+\ldots+h_k=h$ and 
\be
1\leq i_1+h_1<i_2\leq \ldots \ldots\leq i_{k-2}+h_{k-2}< i_{k-1}.\label{eq:pierineq}
\ee
 But now~(\ref{eq:piercorr2})  can be equivalently written as:

\be
\sum_{(h_i,h'')}\ep^{i_1+h_1}\wedge\ldots\wedge\ep^{i_{k-2}+h_{k-2}}\wedge D_{h''}(\ep^{i_{k-1}}\wedge\ep^{i_k}),\label{eq:piercorr3}
\ee
where the sum is over all $(h_1,\ldots,h_{k-2},h")$ such that $h_1+\ldots+h_{k-2}+h''=h$  and satisfying÷
~(\ref{eq:pierineq}).
Since 
\[
D_{h"}(\ep^{i_{k-1}}\wedge\ep^{i_k})=\sum_{\matrix{_{i_{k-1}+h_{k-1}< i_k}\cr _{h_{k-1}+h_{k}=h"}}}\ep^{i_{k-1}+h_{k-1}}\wedge \ep^{i_{k}+h_{k}},
\]
by the inductive hypothesis, substituting into~(\ref{eq:piercorr3}) one gets exactly sum~(\ref{eq:pieruno}).\qed

\medskip
A straightforward application of Pieri's formula~(\ref{eq:pieruno}) gives:
\begin{claim}{\bf Corollary.} \label{corpiergia}
\begin{small}
\[
 D_h(\epsilon^s\wedge\ldots\wedge\epsilon^{s+j-1}\wedge\epsilon^{s+j}\wedge\epsilon^{i_{j+1}}\wedge\ldots\epsilon^{i_k})=\epsilon^s\wedge\ldots\wedge\epsilon^{s+j-1}\wedge D_h(\epsilon^{s+j}\wedge\epsilon^{i_{j+1}}\wedge\ldots\wedge\epsilon^{i_k}). 
\]
\end{small}
\qed
\end{claim}

\claim{}Let $M_n$ be the submodule of $M$ generated by ${\cal E}_n:=(\epsilon^1,\ldots,\epsilon^n)$, $q$  an indeterminate over $\ZZ$ and  $M_n[q]:=M_n\otimes_\ZZ\ZZ[q]$ -- the free $\ZZ[q]$-module spanned by ${\cal E}_n$. As a $\ZZ$-module, the latter is isomorphic to $M$ via the isomorphism 
\[
\left\{\matrix{{\cal Q}_n:&M&\lra &M_n[q]\cr {}&\epsilon^{\alpha\cdot n+i}&\longmapsto&q^\alpha\epsilon^i}\right. ,\quad (\forall\alpha\geq 0,\quad  1\leq i\leq n-1).
\] 
 Let $\wkM_n$ and $\wkM_n[q]\cong \wkM_n\otimes_{\ZZ}\ZZ[q]$ be the $k$-th exterior power of $M_n$ and $M_n[q]$ (thought as a $\ZZ[q]$-module) respectively. Both are freely generated, over $\ZZ$ and $\ZZ[q]$ respectively, by $\{(\ikformep):1\leq i_1<\ldots<i_n\leq n\}$. Let $p_n:\wkM\lra \wkM_n$ be the natural projection defined as:
\[
p_n\left(\sum_{1\leq i_1<\ldots<i_k}a_{i_1\ldots i_k}\cdot\ikformep\right)=\sum_{1\leq i_1<\ldots<i_k\leq n}a_{i_1\ldots i_k}\cdot \ikformep
\]
and $\wedge^k{\cal Q}_n:\wkM\lra\wkM_n[q]$ be the $\ZZ$-module isomorphism induced by ${\cal Q}_n$.
It is easy to see that $p_n\circ D_{h}:\wkM\lra \wkM_n$ is the null homomorphism for all $h\geq n+1$. The proposition below rules the case $h\leq n$. 

\begin{claim}{\bf Corollary.} 
{\em Let $I:=(1\leq i_1<i_2\ldots<i_k\leq n)$ and $0\leq h\leq n$. Then:
\be
p_n\circ D_h(\ikformep)=\sum_{\{{(h_i)\in H(I,h)\,|\, i_k+h_k\leq n}\}} \ihikformep,\label{eq:pierdue}
\ee

and
\begin{eqnarray}
&{}&\wedge^k{\cal Q}_n\circ D_h(\ikformep)=p_nD_h(\ikformep)+\nonumber\\
&{}&+(-1)^{k-1}q\cdot\sum_{\matrix{_{(h_i)\in H(I,h)}\cr _{i_k+h_k-n<i_1}}}\epsilon^{i_k+h_k-n}\wedge\ep^{i_1+h_1}\wedge\ldots\wedge\epsilon^{i_{k-1}+h_{k-1}}.\label{eq:piertre}
\end{eqnarray}
where $H(I,h)$ is as in Theorem~(\ref{Pieri}).
}
\end{claim}
\proof
Equation~(\ref{eq:pierdue}) is obvious: one writes down  expansion~(\ref{eq:pieruno}) and then  projects via $p_n$, canceling all the terms such that $i_k>n$. As for~(\ref{eq:piertre}), one first uses~(\ref{eq:pieruno}) to expand $D_h(\ikformep)$  and then splits the sum as:
\[
D_h(\ikformep)=\sum_{\matrix{_{(h_i)\in H(I,h)}\cr_{i_k+h_k\leq n}}}\ihikformep+\sum_{\matrix{_{(h_i)\in H(I,h)}\cr_{i_k+h_k>n}}}\ihikformep.
\]

\noindent
The first summand occurring on the r.h.s. is precisely $p_nD_h(\ikformep)$.
Applying $\wedge^k{\cal Q}$ to both sides:
\begin{eqnarray}
&{}&\wedge^k{\cal Q}(D_h(\ikformep))=\nonumber\\&{}&=p_nD_h(\ikformep)+ 
\sum_{(h_i)\in H(I,h)}\epsilon^{i_1+h_1}\wedge\ldots\wedge\epsilon^{i_{k-1}+h_{k-1}}\wedge q\epsilon^{i_k+h_k-n}.\,\,\,\,\,\,\,\,\,\,\,\,\,.\label{eq:form6}
\end{eqnarray}

\noindent
Using the $\ZZ_2$-symmetry of $\wedge$, last term of~(\ref{eq:form6}) can be written as \linebreak 
$(-1)^{k-1}q(C+\overline{C})$, where:
\[
(-1)^{k-1}qC:=(-1)^{k-1}q\sum_{\matrix{_{(h_i)\in H(I,h)}\cr _{i_k+h_k-n< i_1}}}\epsilon^{i_k+h_k-n}\wedge\epsilon^{i_1+h_1}\wedge\epsilon^{i_2+h_2}\wedge\ldots\wedge\epsilon^{i_{k-1}+h_{k-1}}\]
is exactly the second summand of the r.h.s. of formula~(\ref{eq:piertre}), while:
\begin{eqnarray}
&{}&\overline{C}:=\sum_{\matrix{_{(h_i)\in H(I,h)}\cr _{i_k+h_k-n\geq i_1}}}\epsilon^{i_k+h_k-n}\wedge\epsilon^{i_1+h_1}\wedge\epsilon^{i_2+h_2}\wedge\ldots\wedge\epsilon^{i_{k-1}+h_{k-1}}=\nonumber\\
&{}&=\sum_{h'=0}^h\sum_{\,\,\,\,\,\,h_k=i_1+n-i_k}^{h'}\ep^{i_k+h_k-n}\wedge\ep^{i_1+h'-h_k}\wedge D_{h-h'}(\ep^{i_2}\wedge\ldots\wedge\ep^{i_k-1}).\,\,\,\,\,\,\,\,\,
\label{eq:form5}
\end{eqnarray}
For each $0\leq h'\leq h$, let $\rho_{h'}$ be the bijection of the set 
\[
\{a\in\NN: i_1+n-i_k\leq a \leq h'\}
\] 
onto itself, defined by $\rho_{h'}(a)=i_1+n+h'-i_k-a$.
Then expression~(\ref{eq:form5}) can also be written as:
\begin{small}
 \begin{eqnarray*}
&{}&\overline{C}=\sum_{h'=0}^h\sum_{\,\,\,\,\,\,h_k=i_1+n-i_k}^{h'}\ep^{i_k+\rho_{h'}(h_k)-n}\wedge\ep^{i_1+h'-\rho_{h'}(h_k)}\wedge D_{h-h'}(\ep^{i_2}\wedge\ldots\wedge\ep^{i_k-1})=\\
&=&\sum_{h'=0}^h\sum_{\,\,\,\,\,\,h_k=i_1+n-i_k}^{h'}\ep^{i_1+h_1}\wedge\ep^{i_k+h_k-n}\wedge D_{h-h'}(\ep^{i_2}\wedge\ldots\wedge\ep^{i_k-1})=-\overline{C}.
\end{eqnarray*}
\end{small}
Thus $\overline{C}=0$ and the proof of~(\ref{eq:piertre}) is complete.
\qed

\claim{}\label{concrem} If one associates to any $\rikformep$ the partition $\ulamb=(r_k,\ldots,r_1)$, then Pieri's formula~(\ref{eq:pierdue}) means precisely to add to the Young diagram $Y(\ulamb)$ of $\ulamb$, contained in a $k(n-k)$ rectangle,  $h$ boxes in all possible ways, no two on the same column (Cf.~(\cite{Fu1}), p.~264): this is a combinatorial way to express classical Pieri's formula holding in the grassmannian $G_k(\CC^n)$ (see also \cite{GH}). Moreover, up to renaming $q$ by $(-1)^{k-1}q$, formula~(\ref{eq:piertre}) is nothing else than quantum Pieri's formula found by Bertram~(\cite{Ber1}). Since $H^*(G_k(\CC^n))$  (resp. $QH^*(G_k(\CC^n)$), the cohomology ring (resp. the small quantum cohomology ring) of $G_k(\CC^n)$, is generated as $\ZZ$-algebra (resp. as $\ZZ[q]$-algebra) by the special Schubert cycles $\sigma_i$ and the product structure is completely determined by Pieri's formula (resp. quantum Pieri's formula), one has hence proven that: 
\begin{claim}{\bf Theorem.}\label{thm3.2} {\em
The cohomology  ring of the grassmannian $G_k(\CC^n)$ (resp. the small quantum cohomology ring) can be realized as a commutative ring of linear operators $\ZZ[D]$ of $\wkM_n$ (resp. $\ZZ[q][D]$ of $\wkM[q]$) via the map $\sigma_i\mapsto D_i$ (resp. $\sigma_i\mapsto D_i$ and $q\mapsto (-1)^{k-1}q$).} \qed
\end{claim}
It is worth to remark that the cohomology rings of $G_k(\CC^n)$, for all  $0\leq k\leq n $, are quotients of the same ring $\ZZ[D]:=\ZZ[D_1,D_2,\ldots]$ of derivations of the exterior algebra $\wM$ of the infinite free $\ZZ$-module $M$.
Once one is given of Pieri's formula and of the Chow basis theorem, everything follows formally (see
e.g.~\cite{Manivel}).  In particular, within our formalism, Giambelli's formula can be recasted as:
\be
\rikformep=\Delta_{(r_k,\ldots,r_1)}(D)\cdot\kformep\quad \forall(r_k\geq\ldots\geq r_1\geq 0)
\ee
where
\[
\Delta_{(r_k,\ldots,r_1)}(D)=\left|\matrix{D_{r_1}&D_{r_2+1}&\ldots&D_{r_k+k-1}\cr
D_{r_1-1}&D_{r_2}&\ldots&D_{r_k+k-2}\cr
\vdots&\vdots&\ddots&\vdots\cr
D_{r_1-k+1}&D_{r_2-k+2}&\ldots&D_{r_k}}\right|,
\]
\noindent setting $D_i=0$ if $i<0$.
 Given any $\ikformep\in \wkM$, Giambelli's problem thus consists in finding $G_{i_1\ldots i_k}(D)\in \ZZ[D]$ (a polynomial expression in $(D_1,D_2,\ldots)$), such that:
\[
\ikformep=G_{i_1\ldots i_k}(D)\cdot\kformep.
\] 
Such a polynomial can be found ``by hands" via suitable  ``integration by parts" (see~\cite{Gat} for details), as indicated in the following simple:
\claim{\bf Example.} Consider $\epsilon^2\wedge\epsilon^5\in\bigwedge^2M$. One has:
\[
\epsilon^2\wedge\epsilon^5=D_1(\epsilon^1\wedge\epsilon^5)-\epsilon^1\wedge\epsilon^6=D_1D_3(\epsilon^1\wedge\ep^2)-D_4(\ep^1\wedge\ep^2)=(D_1D_3-D_4)(\ep^1\wedge\ep^2),
\]
having applied twice Corollary~\ref{corpiergia}.


\begin{thebibliography}{99}



\bibitem{Ber1} A.~Bertram, {\sl Quantum Schubert Calculus},
Adv. Math. {\bf 128}, (1997) 289--305.












\bibitem{Fu1} W.~Fulton,  {\sl Intersection Theory},
Springer-Verlag, 1984.










\bibitem{Gat} L.~Gatto, {\sl Schubert Calculus: An Algebraic Introduction}, 25$^\circ$ Col\'oquio Brasileiro de Matem\'atica, Instituto de Matem\'atica Pura e Aplicada, Rio de Janeiro, 2005.











\bibitem{GH} Ph.~Griffiths, J.~Harris, {\sl Principles of Algebraic
Geometry}, Wiley-Interscience, New York, 1978.








\bibitem{KlLa} S.~L.~Kleiman, D.~Laksov, {\sl Schubert
Calculus}, Amer. Math. Monthly {\bf 79}, (1972), 1061--1082.
















\bibitem{LakTh} D.~Laksov, A.~Thorup, {\sl A Determinantal Formula for the Exterior Powers of the Polynomial Ring}, Preprint, 2004.


\bibitem{LakTh2} D.~Laksov, A.~Thorup, {\sl Private Communication} (forthcoming paper).

\bibitem{Manivel} L.~Manivel, {\sl Fonctions sim\'etriques, polyn\^omes de Schubert et lieux de d\'eg\'en\'erescence}, Cours Sp\'ecialis\'es, SMF, Num\'ero {\bf 3}, 1998.

\bibitem{Mats} H.~Matsumura, {\sl Commutative Rings Theory},
{\bf 8}, Cambridge Univ. Press, Cambridge, 1996.
















\end{thebibliography}
\end{document}